\documentclass[12pt]{article}
\usepackage{amsmath,amsfonts,amssymb,amsthm,color}
\usepackage[dvipdfmx]{graphicx}
\usepackage{bm}
\usepackage{multirow}
\usepackage[margin=2cm]{geometry}
\usepackage{url}

\theoremstyle{plain}
\newtheorem{theorem}{Theorem}[section]

\newtheorem{proposition}[theorem]{Proposition}
\theoremstyle{definition}

\theoremstyle{remark}

\newcommand{\diag}{\mathop{\rm diag}}

\def\pd#1{ \partial_{#1} }
\def\x{\lambda}
\def\y{\tau}

\begin{document}

\title{Holonomic gradient method for distribution function of a weighted sum of noncentral chi-square random variables}

\author{
Tamio Koyama\thanks{Graduate School of Information Science and Technology, University of Tokyo}\ \thanks{Research Fellow of Japan Society for the Promotion of Science}\ \ 
and Akimichi Takemura\footnotemark[1]
}
\date{August, 2015}
\maketitle

\begin{abstract}
We apply the holonomic gradient method to compute the distribution function of 
a weighted sum of independent noncentral chi-square random variables.
It is the distribution function of the squared length
of a multivariate normal random vector.  We treat this distribution as an integral of the
normalizing constant of the Fisher-Bingham distribution on the unit sphere and make use of the
partial differential equations for the Fisher-Bingham distribution.
\end{abstract}

\noindent
{\it Keywords and phrases:} \ 
algebraic statistics,
cumulative chi-square distribution,
Fisher-Bingham distribution,
goodness of fit

\section{Introduction}
\label{sec:intro}
The weighted sum of independent chi-square variables appears in many important problems in statistics.
In the problems for testing against ordered alternatives, cumulative chi-square statistic 
(cf.\ \cite{hirotsu1986}, \cite{nair}) has a good power. For studying the power function of the cumulative chi-square statistic,
we need to evaluate the distribution function of a sum of weighted independent {\it noncentral} chi-square variables.
Goodness of fit test statistics based on empirical cumulative distribution function, such as the
Cram\'er-von Mises statistic or the Anderson-Darling statistic (\cite{anderson-darling}), are infinite sums of weighted independent chi-square variables.
Chapter 4 of \cite{dagostino-stephens} gives a survey of these statistics.  Under an alternative  hypothesis  
the chi-square variables are noncentral.   For studying the power function of these statistics
we want to approximate the infinite sum by a finite sum of sufficiently many terms and compute the cumulative 
distribution of the finite sum.

An exact evaluation of the cumulative distribution function of 
a weighted sum of independent noncentral chi-square random variables was considered to be a difficult numerical problem (see \cite{martinez-blazquez}).  Although the moment generating function is explicitly given, its
Fourier inversion to evaluate the density function and the cumulative distribution function is 
difficult as extensively discussed in Chapter 6 of \cite{tanaka-katsuto}.  See 
\cite{cohen-laplace-book} for the similar problems in other areas of applied mathematics.

Recently in \cite{hgd} we proposed the holonomic gradient method (HGM) for calculating distribution functions and the maximum likelihood estimates using differential equations satisfied by a probability density function with
respect to the parameters.  Since then the method has been successfully used in many problems, including
the computations related to the Fisher-Bingham distribution on the unit sphere 
(\cite{koyama-etal-2014jpaa}, \cite{koyama-2014fe}, \cite{koyama-etal-2014cs}, \cite{sei-kume}).
In this paper we utilize the results on HGM for the Fisher-Bingham distribution
to evaluate the distribution function of a weighted sum of noncentral chi-square random variables.

Let $\mathbf X$ denote a $d$-dimensional random vector following the multivariate normal distribution $N(\mathbf\mu,\Sigma)$. 
Consider the cumulative distribution function $G(r)$ of $\Vert \mathbf X \Vert$: 
\begin{equation}
\label{eq:cdf1}
G(r)=\int_{x_1^2+\cdots +x_d^2 \leq r^2}
\frac{1}{(2\pi)^{d/2}|\Sigma|^{1/2}}
\exp\left(
  -\frac{1}{2}(\mathbf x-\mathbf\mu)^\top\Sigma^{-1}(\mathbf x-\mathbf\mu)
\right)d\mathbf x.
\end{equation}
We call $G(r)$ the ball probability with radius $r$.
By rotation we can assume that $\Sigma=\diag(\sigma_1^2, \dots, \sigma_d^2)$ is a diagonal matrix without loss of generality.  Hence
$G(r)$ is the distribution function of the square root of a  weighted sum of independent 
noncentral chi-square random variables, where weights are $\sigma_i^2$, $i=1,\dots,d$.
Furthermore the conditional distribution of $\mathbf X$ given its length $r=\Vert \mathbf X \Vert$
is the Fisher-Bingham distribution. 
This fact allows us to directly apply the results
for the Fisher-Bingham distribution to the evaluation of the distribution of the weighted sum of independent noncentral chi-square random variables.  
As we show in Section \ref{sec:numerical} our method works very well, both in accuracy and speed.

The organization of this paper is as follows.  
In Section \ref{sec:hgm}  we summarize known results on  HGM
for the Fisher-Bingham distribution and show how they can be used to evaluate
the distribution of the a  weighted sum of independent noncentral chi-square random variables. We also discuss the problem of initial values needed to use HGM.
In Section \ref{sec:laplace} we present asymptotic results for the Fisher-Bingham integral and its derivatives for the
case that the length of the multivariate normal vector diverges to infinity.  This result is used to check the the numerical  accuracy of our experiments in Section \ref{sec:numerical}.
We end the paper with some discussions in Section \ref{sec:summary}.

\medskip\noindent
{\bf Acknowledgment}\/.
This work is supported by
JSPS Grant-in-Aid for Scientific Research No.\ 25220001
and Grant-in-Aid for JSPS Fellows No.\ 02603125.

\section{Holonomic system and initial values}
\label{sec:hgm}
Let 
\begin{align*}
&\Sigma=\diag(\sigma_1^2, \dots, \sigma_d^2), & \quad \mathbf\mu=(\mu_1,\dots,\mu_d)^\top. \\
\end{align*}
We define new parameters $\lambda_i, \tau_i$, $i=1,\dots,d$, by
\begin{align*}
&\x_i=-\frac{1}{2\sigma_i^2}, & 
\y_i=\frac{\mu_i}{\sigma_i^2}
\end{align*}
and the Fisher-Bingham integral $f(\x,\y,r)$ by
\begin{equation}
f(\x,\y,r)
=\int_{S^{d-1}(r)} \exp\left(\sum_{i=1}^d\x_it_i^2+\sum_{i=1}^d\y_it_i \right)d\mathbf t, 
\label{eq:fb-integral}
\end{equation}
where $\lambda=(\lambda_1,\dots,\lambda_d)$, $\tau=(\tau_1, \dots, \tau_d)$, 
$
S^{d-1}(r)
=\left\{\mathbf t\in\mathbf R^d \mid t_1^2+\cdots +t_d^2=r^2\right\}
$ 
is the sphere of radius $r$ and
$d\mathbf t$ is the volume element of $S^{d-1}(r)$ so that
\[
\int_{S^{d-1}(r)} d\mathbf t = r^{d-1} S_{d-1}, \quad S_{d-1}= {\rm Vol}(S^{d-1}(1))=\frac{2 \pi^{d/2}}{\Gamma(d/2)}.
\]
Then $G(r)$ in  \eqref{eq:cdf1} is written as
\begin{equation}
\label{eq:cdf2}
G(r)=
\frac{\prod_{i=1}^d\sqrt{-\x_i}}{\pi^{d/2}}
\exp\left(
  \frac{1}{4}\sum_{i=1}^d \frac{\y_i^2}{\x_i}
\right)
\int_0^r f(\x,\y,s)ds.
\end{equation}
We will numerically integrate the right-hand side of \eqref{eq:cdf2}.
We denote the partial differential operator with respect to $\x$ by $\pd{\x}$.
For $\mathbf t\in S^{d-1}(r)$,  $(t_1^2 + \dots + t_d^2)/r^2=1$ and 
\begin{align}
f(\x,\y,r) 
&= 
\int_{S^{d-1}(r)} \frac{1}{r^2} (t_1^2 + \dots + t_d^2) \exp\left(\sum_{i=1}^d\x_it_i^2+\sum_{i=1}^d\y_it_i \right)d\mathbf t
\nonumber \\
&=
\frac{1}{r^2}
\left(  \pd{\x_1}+\cdots+\pd{\x_d}\right) f(\x,\y,r).
\label{eq:sum-ys}
\end{align}
By HGM we evaluate $\pd{\x_i}f(\x,\y,r)$, $i=1,\dots,d$, and use
\eqref{eq:sum-ys} to compute $f(\x,\y,r)$. 
In fact we also evaluate $\pd{\y_i}f(\x,\y,r)$, $i=1,\dots,d$.

Define a $2d$-dimensional vector of partial derivatives of $f(\x,\y,r)$ by
\begin{equation}
\label{eq:vector-F}
\mathbf F = \left(
  \pd{\y_1}f, \dots, \pd{\y_d}f,
  \pd{\x_1}f, \dots,\pd{\x_d}f 
\right)^\top.
\end{equation}
Elements of $\mathbf F$ are called ``standard monomials'' in HGM.
By Theorem 3 of \cite{koyama-etal-2014cs} we have
\begin{equation}
\label{eq:pfaffian-r}
\pd{r}  \mathbf F = P_r \mathbf F,
\end{equation}
where the  $2d \times 2d$  matrix $P_r= (p_{ij})$, called the Pfaffian matrix, is of the form
\begin{equation}
\label{eq:pfaff-1}
P_r = 
\frac{1}{r} 
\left(
\begin{array}{cccccc}
2r^2\x_1+1 & & {\bm O} & \y_1 & \cdots & \y_1 \\
&\ddots & & &\vdots  & \\
{\bm O} & &2r^2\x_d+1 &  \y_d &\cdots & \y_d \\
   r^2\y_1 &    &    {\bm O} & 2r^2\x_1+2 &   &      {\bm 1} \\
&\ddots & & &\ddots &\\
      {\bm O} &&    r^2\y_d &        {\bm 1}& & 2r^2\x_d+2 \\  
\end{array}
\right), 
\end{equation}
with ${\bm O}$ denoting an off-diagonal block of 0's and ${\bm 1}$ denoting an off-diagonal block of 1's.
The elements $p_{ij}$ of $P_r$ are expressed as 
\begin{align*}
rp_{ij}
&= (2\x_ir^2+1)\delta_{ij} + \sum_{k=1}^{d} \y_i\delta_{j(k+d)}
\quad (1 \leq i \leq d ),\\
rp_{(i+d)j} 
&= \y_ir^2\delta_{ij} + (2\x_ir^2+2)\delta_{j(i+d)}+ \sum_{k\neq i}\delta_{j(k+d)}
 \quad (1 \leq i \leq d ), 
\end{align*}
for $1 \leq j \leq 2d$, where $\delta_{ij}$ denotes Kronecker's delta.
Given initial values for the elements of $\mathbf F$ at $r=r_0$, we can apply a standard ODE solver 
to \eqref{eq:pfaffian-r} for numerically evaluating  $\mathbf F$.

For the initial values at a small $r=r_0>0$, we can use the following series expansion
of the Fisher-Bingham integral (\cite{Kume-Walker}):
\begin{equation}  \label{eq:series}
f(\x,\y,r)
=
r^{d-1}S_{d-1} \times
\sum_{\alpha, \beta \in {\mathbb N}_0^d}
r^{2|\alpha+\beta |}
\frac{(d-2)!!\prod_{i=1}^{d}(2\alpha_i+2\beta_i-1)!!}{(d-2+2|\alpha|+2|\beta|)!!\alpha !(2\beta)!}
\x^\alpha \y^{2\beta},
\end{equation}
where ${\mathbb N}_0 = \{0, 1, 2, \dots \}$ and 
for a multi-index $\alpha \in {\mathbb  N}_0^d$ we define 
\[
\alpha! = \prod_{i=1}^d\alpha_i!,\quad  \alpha!! = \prod_{i=1}^d\alpha_i!! \ 
\text{and} \ |\alpha| = \sum_{i=1}^d\alpha_i.
\]
By term by term differentiation of this series we can evaluate derivatives of $f(\x,\y,r)$.
For computing the initial values, we apply the following approximation:
\begin{align}
\label{initial1}
\frac{\partial f}{\partial \y_i} &= S_{d-1}r^{d+1} \y_i + O(r^{d+3})
\quad (i=1,\dots, d),\\
\label{initial2}
\frac{\partial f}{\partial \x_i} &= S_{d-1}r^{d+1}  + O(r^{d+3})
\quad (i=1,\dots, d).
\end{align}
By this approximation, we reduce the computational time for the initial values.
However the accuracy of the result does not decrease at all 
as we will show in Section \ref{sec:numerical}.

As $r\rightarrow\infty$, the absolute values of $f(\x,\y,r)$ and its derivatives become
exponentially small, as we analyze the behavior in the next section.
Hence we also consider the following vector
\begin{equation}
\label{eq:large-r-1}
\mathbf Q=
\exp(-r^2\x_1-r|\y_1|)
\left(
\frac{1}{r}\pd{\y_1}f, 
\pd{\y_2}f, \dots, \pd{\y_d}f,
\frac{1}{r^2}\pd{\x_1}f, 
\pd{\x_2}f, \dots, \pd{\x_d}f
\right)^\top.
\end{equation}
Then from \eqref{eq:pfaffian-r} it is easy to obtain $\pd{r} \mathbf Q$  as
\begin{equation}\label{eq:large-r-2}
\pd{r}  \mathbf Q
= \left(D^{-1}\pd{r} D - (2r\x_1+|\y_1|)I_{2d} + DP_rD^{-1} \right)\mathbf Q,
\end{equation}
where 
$I_{2d}$ is the identity matrix with size $2d$ 
and 
$$
D = \diag\left(\frac{1}{r},1,\dots,1,\frac{1}{r^2},1,\dots,1\right).
$$
The equation \eqref{eq:large-r-1} is a refinement of the equation $(21)$ in 
\cite{koyama-etal-2014cs}.
By Proposition \ref{prop:lap1} in the next section, 
each element of $\mathbf Q$ converges to some non-zero value when $r$ goes to 
the infinity.
This prevents the adaptive Runge-Kutta method from slowing down.



\section{Laplace approximation close to the infinity}
\label{sec:laplace}

In our implementation of HGM, we start from a small $r=r_0 > 0$ and numerically integrate $\mathbf F$ 
in \eqref{eq:vector-F} up to $r=1$ and then integrate $\mathbf Q$ in \eqref{eq:large-r-1} toward $r=\infty$.  
In order to assess the accuracy of $\mathbf Q$ for large $r$, we derive the asymptotic values of the elements of $\mathbf Q$ by the Laplace method.
The Laplace approximation, including higher order terms, for the Fisher-Bingham integral itself was given in \cite{kume-wood}. However here we also need approximations for its derivatives, which were not given in
\cite{kume-wood}. Hence we give the approximations of the main terms of the Fisher-Bingham integral
and its derivatives and a sketch of their proofs.

We first consider the case of single largest $\x_1$.  We state the following result.
\begin{proposition}  
\label{prop:lap1}
Suppose $ %
\label{eq:ass1}
0 > \x_1 > \x_2 \ge \dots \ge \x_d.
$
Then, as $r\rightarrow\infty$,
\begin{align}
f(\x,\y,r)&= \frac{\pi^{(d-1)/2}}{\prod_{i=2}^d(\x_1-\x_i)^{1/2} }(e^{r \y_1} + e^{- r \y_1})
\exp\left(r^2 \x_1  -\sum_{i=2}^d \frac{\y_i^2}{4(\x_i-\x_1)}\right) (1+o(1)).
\label{eq:lap-f-single}\\
\pd{\x_1}f(\x,\y,r) &= r^2 f(\x_1,\y_1, r) (1+o(1)), 
\label{eq:lap-fx1-single} \\
\pd{\x_j}f(\x,\y,r) & =\left\{  \left(    \frac{\y_j}{2(\x_j-\x_1)}
  \right)^2  +\frac{1}{2(\x_1-\x_j)} \right\} f(\x,\y,r) (1+o(1)),  \quad (j=2,\dots, d)
\label{eq:lap-fxi-single}
\\
\pd{\y_1}f(\x,\y,r)& =   r\frac{e^{r\y_1}-e^{-r\y_1}}{e^{r\y_1}+e^{-r\y_1}} f(\x,\y,r) (1+o(1)),
\label{eq:lap-fy1-single} \\
\pd{\y_j}f(\x,\y,r)& =   \frac{\y_j}{2(\x_1-\x_j)} f(\x,\y,r)(1+o(1)), \quad (j=2,\dots, d).
\label{eq:lap-fyi-single} 
\end{align}
\end{proposition}

Note that for $\y_1 > 0$, in  
\eqref{eq:lap-f-single} $e^{-r\y_1}$ is exponentially smaller than
$e^{r\y_1}$ and it can be omitted.  However we leave $e^{-r\y_1}$ there for consistency with the case
of $\y_1=0$. Also we found that leaving $e^{-r\y_1}$ in \eqref{eq:lap-f-single} greatly improves the approximation.

We now give a rough proof of Proposition \ref{prop:lap1}.  In the proof, the main contributions
from the neighborhoods of maximal points are carefully evaluated, but the contributions from
outside the neighborhoods are not bounded rigorously.
Replacing $t_i$ by  $r t_i$ and
integrating over $S^{d-1}(1)$ can write
\begin{align}
f(\x,\y,r)&=r^{d-1} \int_{S^{d-1}(1)} \exp\left(r^2\sum_{i=1}^d \x_it_i^2+r \sum_{i=1}^d\y_it_i \right)d\mathbf t,
\label{fb1} \\
\partial_{\x_j} f(\x,\y,r)&=r^{d+1} \int_{S^{d-1}(1)} t_j^2\exp\left(r^2\sum_{i=1}^d \x_it_i^2+r\sum_{i=1}^d \y_it_i \right)d\mathbf t,
\label{fb2} \\
\partial_{\y_j} f(\x,\y,r)&=r^d \int_{S^{d-1}(1)} t_j\exp\left(r^2\sum_{i=1}^d \x_it_i^2+r\sum_{i=1}^d \y_it_i \right)d\mathbf t .
\label{fb3} 
\end{align}

For very large $r$
\begin{equation}
\label{eq:lap1}
r^2 (\x_1 t_1^2 + \x_2 t_2^2 + \dots + \x_d t_d^2),  \qquad 1=t_1^2 + \dots + t_d^2,
\end{equation}
takes its maximum value at two points $t_1=\pm 1, t_2=\dots=t_d=0$.
The main contributions to \eqref{fb1}--\eqref{fb3} 
come from neighborhoods of these two points
$(\pm 1, 0,\dots,0)$. The contribution from the complement of these two neighborhoods should be 
exponentially small as $r\rightarrow\infty$, although we do not give a detailed argument.
We also have to consider the effect of
$r \sum_{i=1}^d \y_i t_i$.   But it is of the order $O(r)$, whereas \eqref{eq:lap1}
is of the order $O(r^2)$.  Hence $r \sum_{i=1}^d \y_i t_i$ only perturbs the maximizing values
$(\pm 1, 0,\dots,0)$ by the term of the order $O(1/r)$.  Based on these considerations write
\[
t_1^2 = 1-t_2^2 - \dots -t_d^2,   \qquad t_1 = \pm \sqrt{1-t_2^2 - \dots -t_d^2} 
\doteq \pm \big(1 - \frac{1}{2}(t_2^2 + \dots + t_d^2)), 
\]
where $|t_2|, \dots, |t_d|$ are small.  As shown below, $|t_i|$, $i=2,\dots,d$, are of the order $O(1/r)$.
We now consider the neighborhood of  $(1,0,\dots,0)$. By completing the squares we have
\begin{align}
&r^2 \sum_{i=1}^d \x_it_i^2+ r \sum_{i=1}^d \y_i t_i \nonumber \\
&\qquad = r^2 \x_1 + r \y_1 +  r^2 \sum_{i=2}^d \big((\x_i - \x_1 -\frac{\y_1}{2r}) t_i^2 + \frac{\y_i}{r} t_i\big)  + o(1)
\nonumber \\&\qquad 
=r^2 \x_1 + r \y_1 + \sum_{i=2}^d 
\big[ (\x_i - \x_1 - \frac{\y_1}{2r}) (rt_i + \frac{\y_i}{2(\x_i - \x_1- \frac{\y_1}{2r})})^2
- \frac{\y_j^2}{4(\x_i - \x_1 -\frac{\y_1}{2r})}\big] + o(1) 
\label{eq:asymp-exp1}
\\&\qquad 
= r^2 \x_1 + r \y_1 + \sum_{i=2}^d 
\big[ (\x_i - \x_1)(rt_i + \frac{\y_i}{2(\x_i - \x_1)})^2
- \frac{\y_j^2}{4(\x_i - \x_1)}\big] + o(1) .\nonumber
\end{align}
Furthermore around $(1,0,\dots,0)$ the volume element $d\mathbf t$ of the unit sphere $S^{d-1}(1)$ is 
approximately equal to the Lebesgue measure $dt_2 \dots dt_d$,  with the error of 
the order $t_2^2 + \dots + t_d^2$.  Hence by the change of variables
\[
u_i = r t_i, \quad i=2,\dots,d,
\]
the contribution to $f(\x,\y,r)$ from the neighborhood of $(1,0,\dots,0)$ is evaluated as
\begin{align}
&\exp(r^2 \x_1 + r \y_1) \int_{{\mathbb R}^{d-1}} \exp( (\x_i - \x_1)(u_i + \frac{\y_i}{2(\x_i - \x_1)})^2
- \frac{\y_j^2}{4(\x_i - \x_1)}) du_1 \dots du_d \nonumber \\
& \qquad
= 
\exp\left(  r^2   \x_1 +  r \y_1  -\sum_{i=2}^d \frac{\y_i^2}{4(\x_i-\x_1)}
\right)
\frac{\pi^{(d-1)/2}}{\prod_{i=2}^d(\x_1-\x_i)^{1/2} }.
\label{eq:single-asym-1}
\end{align}

Similarly by changing the sign of $\tau_1$ we can evaluate 
the contribution from the neighborhood of $(-1,0,\dots,0)$ as
\begin{equation}
\exp\left(  r^2   \x_1 -  r \y_1  -\sum_{i=2}^d \frac{\y_i^2}{4(\x_i-\x_1)}
\right)
\frac{\pi^{(d-1)/2}}{\prod_{i=2}^d(\x_1-\x_i)^{1/2} }.
\label{eq:single-asym-2}
\end{equation}
Adding  \eqref{eq:single-asym-1} and \eqref{eq:single-asym-2} we obtain
\eqref{eq:lap-f-single}.  

For $\pd{\x_1} f(\x,\y,r)$ and 
$\pd{\y_1} f(\x,\y,r)$, we can just put $t_1 = \pm 1$ in
\eqref{fb2} and \eqref{fb3}.  Adding contributions from 
two neighborhoods we obtain \eqref{eq:lap-fx1-single} and \eqref{eq:lap-fy1-single}.

For $\pd{x_i} f(\x,\y,r)$ and $\pd{\y_i} f(\x,\y,r)$,
$j\ge 2$, 
we write
\begin{align*}
t_j &= \frac{u_j}{r} = \frac{1}{r} \left(u_j + \frac{\y_j}{2(\x_j - \x_1)} - \frac{\y_j}{2(\x_j - \x_1)}\right), \\
t_j^2 &= \frac{1}{r^2} \left(u_j + \frac{\y_j}{2(\x_j - \x_1)} - \frac{\y_j}{2(\x_j - \x_1)}\right)^2
\end{align*}
and take the expectation with respect to a normal density.  Then we obtain
\eqref{eq:lap-fxi-single} and \eqref{eq:lap-fyi-single}.  Although we did not
give a detailed analysis of the remainder terms, we can show that the relative errors
in \eqref{eq:lap-fx1-single}--\eqref{eq:lap-fyi-single} are of the order $O(1/r)$.
This completes the proof of Proposition \ref{prop:lap1}.

A generalization of Proposition \ref{prop:lap1} to the case that
$\x_1 = \dots = \x_m > \x_{m+1} \ge \dots \ge \x_d$ is given in Appendix.
We note that numerically HGM works fine even if 
some of the $\lambda$'s are close to one another, because  the Pfaffian system does not have a singular locus
except at $r=0$ and the main exponential order is the same
in Proposition \ref{prop:lap1} and in Proposition \ref{prop:lap-2}.
However when we want to check whether the ratio of HGM to the asymptotic value is close to one, 
then we have difficulty when some of the $\lambda$'s are close to one another.

\section{Numerical experiments}
\label{sec:numerical}

In this section we describe our numerical experiments on the performance of HGM.
The programs and the raw data of our numerical experiments are obtained 
at 
\begin{center}
{\tt http://github.com/tkoyama-may10/ball-probability/}
\end{center}
Our programs utilize the Gnu Scientific Library\cite{gsl}.

In our experiments we compute the initial values of $r^{-(d+1)}\mathbf F$ at 
$r=r_0=1.0\times 10^{-6}$ by 
\eqref{initial1} and \eqref{initial2}.
The reason for multiplying $\mathbf F$ by $r^{-(d+1)}$ is 
that the values of elements of $\mathbf F$ are too small at $r_0$ for floating point numbers.
Then up to $r=1$, we solve 
the differential equation \eqref{eq:pfaffian-r} numerically.
In our implementation, we utilize 
explicit embedded Runge-Kutta Prince-Dormand (8, 9) method and we set the accuracy to $1.0\times 10^{-6}$.
In order to prevent the elements of $\mathbf F$ becoming too large, 
we re-scale the elements of $\mathbf F$ several times.
Then at $r=1$ we switch to $\mathbf Q$ in \eqref{eq:large-r-1} and solve \eqref{eq:large-r-2}.

Note that we can not take $r=0$ as an initial point.
The point $r=0$ is in the singular locus of the differential
equation \eqref{eq:pfaffian-r} 
since the denominator of $P_r$ becomes zero at the point.
Hence, numerical differential equation solvers can not compute
the differentiation of $\mathbf F$ at the point.

Our implementation computes the initial value of $\mathbf F$ 
by the approximations \eqref{initial1} and \eqref{initial2}, 
which use only the first term of the series expansions.
Hence, we have to take very small value for $r$ in order to reduce the error.

Each component of vector $\mathbf F$ takes very small value at $r_0=10^{-6}$.
We deal with this problem by storing not the value of $\mathbf F$ itself but 
the product of $\mathbf F$ and a large constant in double precision type array.
Related to this problem, there is another problem that
the value of each component of $\mathbf F$ increases rapidly
when we solve ordinary differential equation \eqref{eq:pfaffian-r}  numerically.
We multiply vector $\mathbf F$ by a small constant when a component of $\mathbf F$ becomes
larger than a fixed value.
By this way, our implementation prevents values in double precision type array 
becoming too large.

\begin{figure}[htbp]
\begin{center}
\includegraphics[width=0.35\hsize]{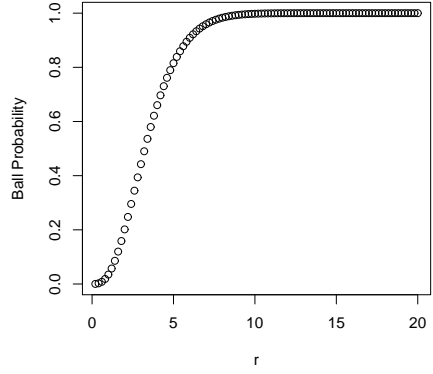}
\caption{ CDF  $G(r)$ for the first experiment}
\label{fig1}
\end{center}
\end{figure}

Our first experiment is for $d=3$ and the following parameter values
\begin{align}
\sigma_1 &= 3.00, & \sigma_2 &= 2.00, & \sigma_3 &= 1.00, \nonumber\\
\mu_1 &= 1.00, & \mu_2 &= 0.50, & \mu_3 &= 0.25, \label{eq:5}
\end{align}
i.e., 
\begin{align*}
\x_1 &= -0.0555556, & \x_2 &= -0.125, & \x_3 &= -0.5, \\
\y_1 &= 0.111111, & \y_2 &= 0.125, & \y_3 &= 0.25.
\end{align*}
By HGM we compute $G(r)$. We show its graph in Figure \ref{fig1}
to confirm that our implementation correctly calculated 
the asymptotic behavior as $G(r)\rightarrow 1$ as $r\rightarrow\infty$.

\begin{figure}[htbp]
\begin{center}
\includegraphics[width=0.45\hsize]{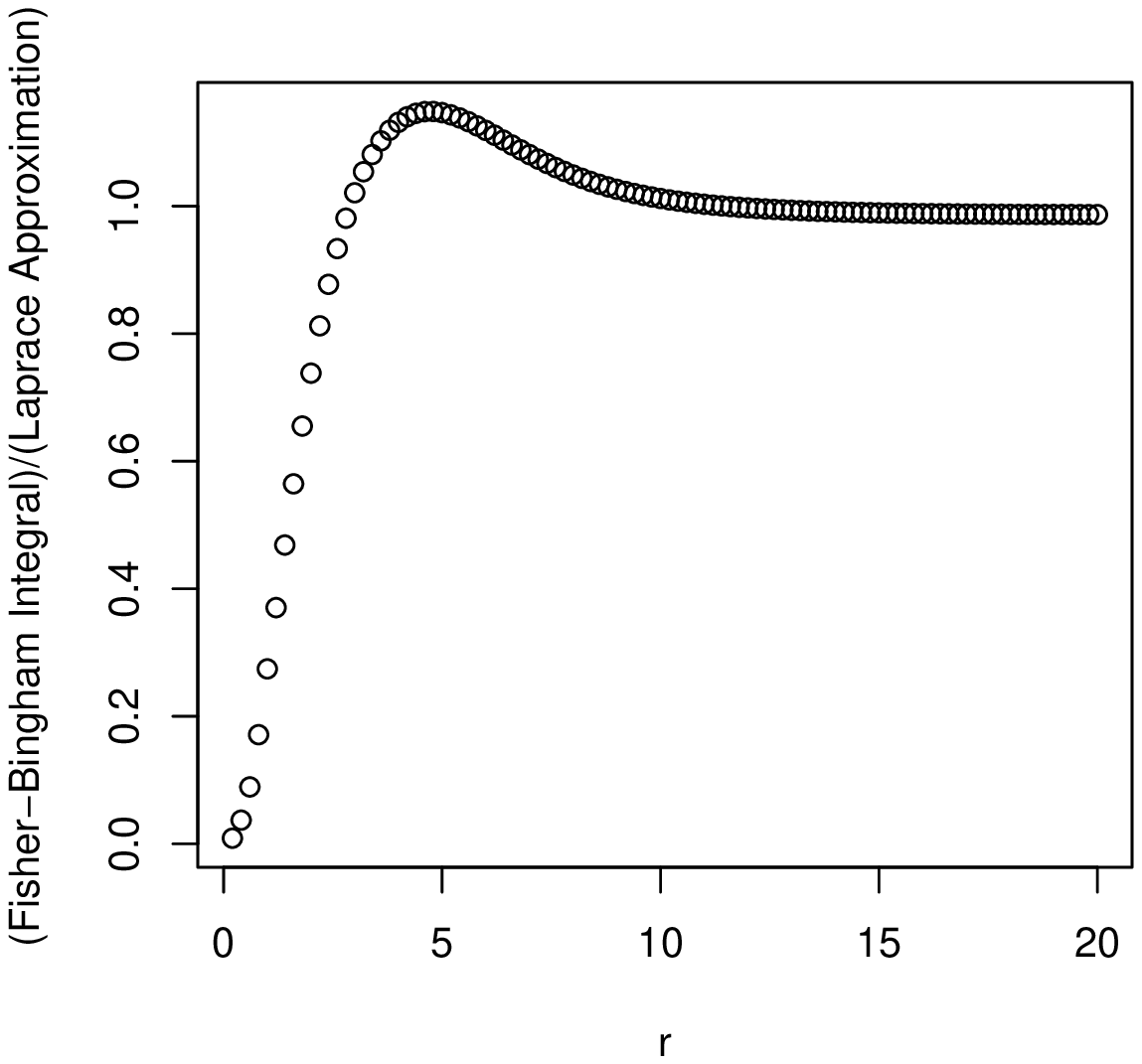} \qquad
\includegraphics[width=0.45\hsize]{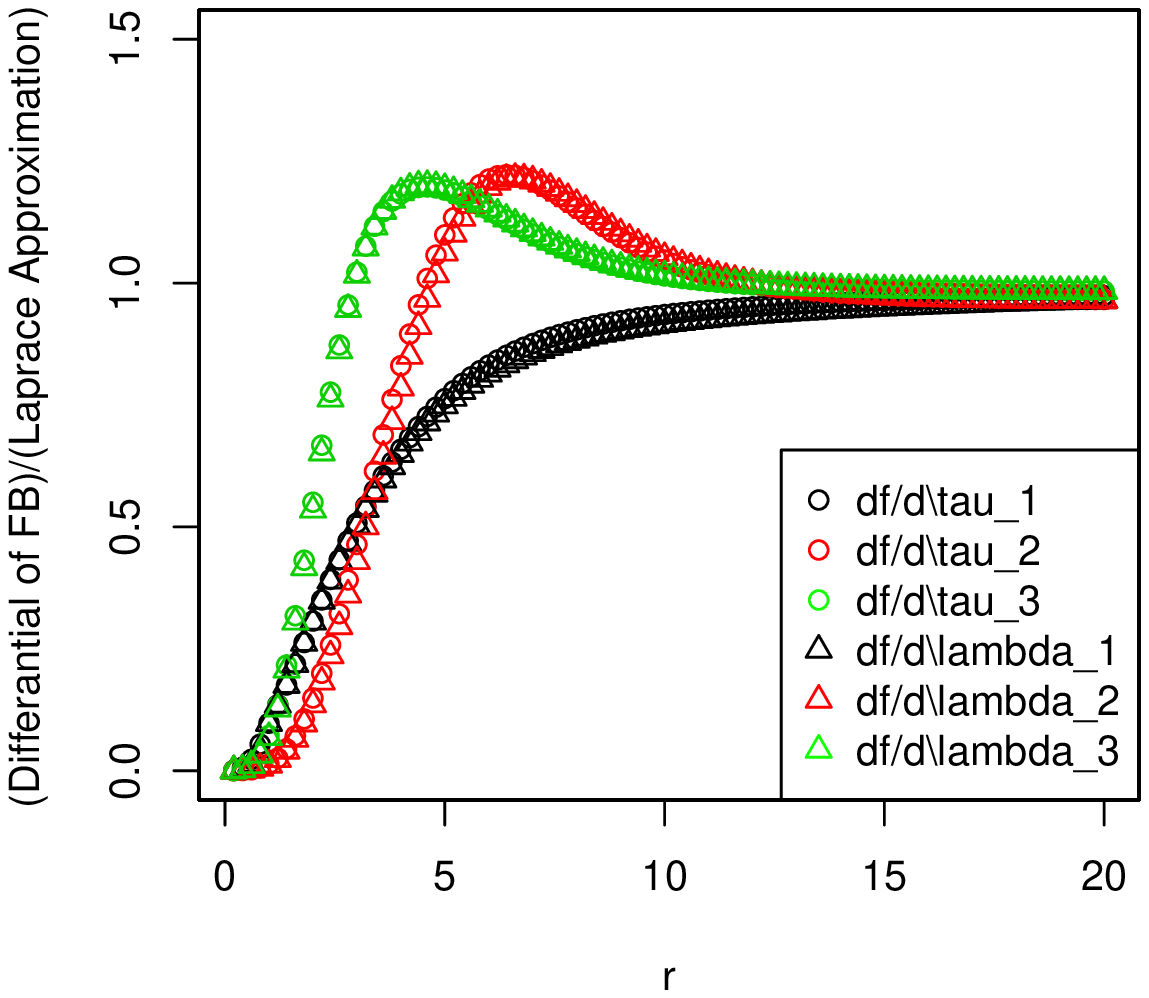}
\caption{Ratios to the Laplace approximations}
\label{fig2}
\end{center}
\end{figure}

For this example, we also check the accuracy by computing the ratios of
$f(\x,\y,r)$ and the elements of $\mathbf F$ to their asymptotic expressions in Proposition
\ref{prop:lap1}. 
The left figure of Figure \ref{fig2} shows the ratio of $f(\x,\y,r)$ 
to its asymptotic expression and the right figure shows the ratios of elements of $\mathbf F$ to their asymptotic expressions.  
Note that the value of the ratio corresponding to $\partial f/\partial\lambda_i$
is very close to that of $\partial f/\partial\tau_i$ so that  
the triangles overlap with  the circles.
We see that the numerical integration involved in HGM, starting from a small $r_0$, 
is remarkably accurate, so that the ratios numerically converge to 1 as $r\rightarrow\infty$.

In our second example we consider diagonal matrices $\Sigma^{(1)}$ and $\Sigma^{(2)}$ 
with diagonal elements
\begin{equation}\label{kym-hirotsu-1}
(\sigma_k^2)^{(1)} = \frac{d+1}{k(k+1)} \quad (1\leq k\leq d),
\end{equation}
and 
\begin{equation}\label{kym-hirotsu-2}
(\sigma_k^2)^{(2)} = \frac{2(d+2)(d+3)}{k(k+1)(k+2)(k+3)} \quad (1\leq k\leq d),
\end{equation}
respectively.  
These weights are considered for cumulative chi-square statistics in \cite{hirotsu1986}.
Let 
\begin{align*}
\mathbf\mu^{(1)} &= 0,\\
\mathbf\mu^{(2)} &= 
\begin{pmatrix}
0 & 0.01 & 0.02 & \cdots & 0.01\times (d-1)
\end{pmatrix}^\top.
\end{align*}

For each dimension $d$,  we computed the probability  
$P\left(10^{-6}\leq \Vert {\mathbf X}\Vert< 40.0\right)$
and measured the computational times in seconds.
We considered the following four patterns of  parameters:
\begin{align*}
&(\Sigma^{(1)},\mathbf\mu^{(1)}),&
&(\Sigma^{(1)},\mathbf\mu^{(2)}),\\
&(\Sigma^{(2)},\mathbf\mu^{(1)}),&
&(\Sigma^{(2)},\mathbf\mu^{(2)}).
\end{align*}
The experimental results are shown in Table \ref{tab1}. 
$1-p$  stands for the values $1-P\left(10^{-6}\leq \Vert {\mathbf X}\Vert< 40.0\right)$
are generally accurate to $10^{-8}$.

\begin{table}[hbp]
\caption{Accuracy and computational times for $\Sigma^{(1)}$ and $\Sigma^{(2)}$} 
\label{tab1}
\smallskip
\centering
\begin{tabular}{c|cccccccc}
\hline
\multirow{3}{*}{dimension}&
\multicolumn{4}{c}{$\Sigma^{(1)}$}&
\multicolumn{4}{c}{$\Sigma^{(2)}$}\\
&
\multicolumn{2}{c}{$\mathbf\mu=0$}&
\multicolumn{2}{c}{$\mathbf\mu\neq 0$}&
\multicolumn{2}{c}{$\mathbf\mu=0$}&
\multicolumn{2}{c}{$\mathbf\mu\neq 0$}\\
& $1-p$ & times(s) & $1-p$ & times(s)&  $1-p$ & times(s) & $1-p$ & times(s) \\\hline
10 & 1.60e-08 & 0.03 & 1.60e-08 & 0.03   & 1.60e-08 & 0.11 & 2.10e-09 & 0.11 \\
11 & 1.76e-08 & 0.03 & 1.57e-08 & 0.04   & 1.76e-08 & 0.12 & 1.56e-09 & 0.14 \\
12 & 1.61e-08 & 0.04 & 1.15e-08 & 0.04   & 1.61e-08 & 0.16 & 9.59e-10 & 0.17 \\
13 & 1.81e-08 & 0.04 & 1.05e-08 & 0.04   & 1.80e-08 & 0.20 & 7.90e-10 & 0.19 \\
14 & 2.02e-08 & 0.04 & 9.95e-09 & 0.05   & 2.02e-08 & 0.24 & 6.94e-10 & 0.25 \\
15 & 2.34e-08 & 0.04 & 9.58e-09 & 0.06   & 2.34e-08 & 0.30 & 6.44e-10 & 0.30 \\
16 & 2.77e-08 & 0.06 & 9.73e-09 & 0.07   & 2.77e-08 & 0.36 & 2.89e-10 & 0.36 \\
17 & 3.40e-08 & 0.07 & 4.85e-09 & 0.08   & 3.40e-08 & 0.41 & 2.74e-10 & 0.42 \\
18 & 1.89e-08 & 0.08 & 4.62e-09 & 0.08   & 1.89e-08 & 0.49 & 2.82e-10 & 0.52 \\
19 & 2.08e-08 & 0.08 & 4.40e-09 & 0.10   & 2.09e-08 & 0.56 & 4.05e-10 & 0.57 \\
20 & 2.33e-08 & 0.10 & 4.32e-09 & 0.11   & 2.41e-08 & 0.65 & 1.13e-09 & 0.65 \\
\hline
\end{tabular}
\end{table}

As the radius $r$ increases or the dimension $d$ of the sphere increases,
our implementation takes long time to evaluate.
Table 1 shows that 
the computational complexity also depends on the values of $\x$.
However we do not know what value of $\x$ makes the computational time worse.

As our third example we consider how our method works for large dimension.
Corresponding to the asymptotic null distribution of 
Anderson-Darling statistic, which is an infinite sum of weighted
$\chi^2$ variables, consider the weights
$$
\sigma_k^2 = \frac{1}{k(k+1)}, \ \  \mu_k = 0 \quad (1\leq k\leq d).
$$
Here we truncate the infinite series at $d$.
We computed the probability and measured its computational time.
We fixed the radius as $r=20.0$. The results on the computational time 
are shown in Table \ref{tab:AD} and its figure. Even for $d=100$, our method is accurate and fast enough to be practical.
This is a remarkable progress since the implementation of HGM in 
\cite{koyama-etal-2014cs} can compute only up to dimension $d=8$.
The key idea for this progress are 
the simple approximation of the initial values \eqref{initial1} and \eqref{initial2} 
for HGM and 
the refined differential equation \eqref{eq:large-r-2} based on 
the Laplace approximation.

The computational bottleneck of HGM is the computation of $P_r{\mathbf F}$ in each
step of solving the ODE. By the form of the matrix $P_r$, the number of additions
in each step increases in order $O(d^2)$. We guess this is a reason 
that growth of computational times in the figure of Table \ref{tab:AD}
seems to be in the order $O(d^2)$.

\begin{table}[htbp]
\caption{Computational times for Anderson-Darling statistic}
\label{tab:AD}
\smallskip
\begin{minipage}[c]{0.45\hsize}
\begin{center}
\begin{tabular}{ccc}
dim& $1-p$ & time(s)\\
\hline
30 & 5.70e-08 & 1.03 \\ 
35 & 3.76e-08 & 1.59 \\ 
40 & 4.85e-08 & 2.36 \\ 
45 & 6.13e-08 & 3.30 \\ 
50 & 8.97e-08 & 4.42 \\ 
55 & 5.29e-08 & 5.94 \\ 
60 & 7.91e-08 & 7.56 \\ 
65 & 6.28e-08 & 9.69 \\ 
70 & 1.02e-07 & 12.05 \\ 
75 & 6.77e-08 & 14.63 \\ 
80 & 7.22e-08 & 17.81 \\ 
85 & 6.25e-08 & 21.33 \\ 
90 & 5.64e-08 & 25.10 \\ 
95 & 5.21e-08 & 29.54 \\ 
100 & 4.90e-08 & 35.05 \\ 
\hline
\end{tabular}
\end{center}
\end{minipage}
\begin{minipage}[c]{0.45\hsize}
\begin{center}
\includegraphics[width=0.95\hsize]{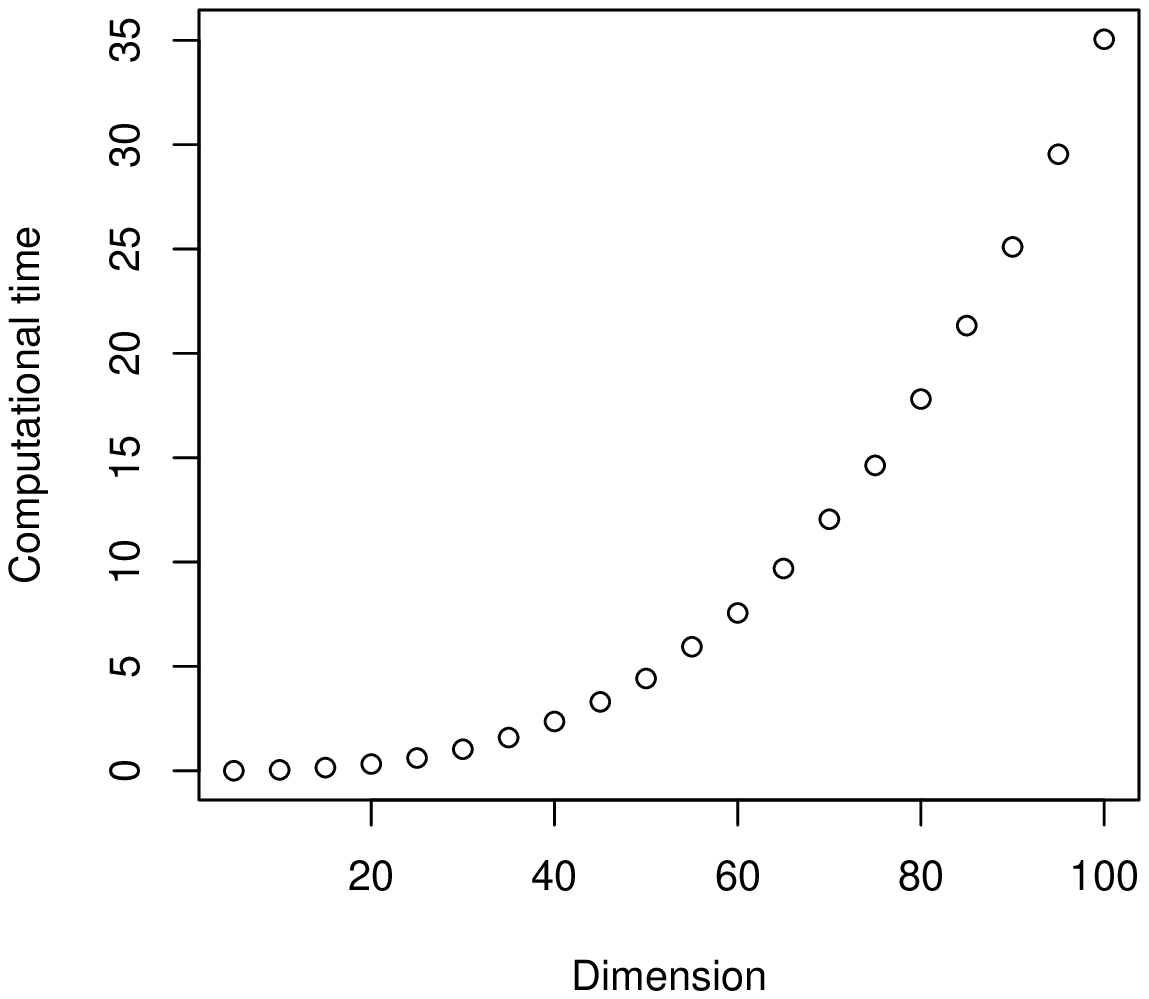}\\
Graph of computational times
\end{center}
\end{minipage}
\end{table}

As our fourth example we consider the case where $\Sigma$ is the identity matrix
and $\mu=0$. 
In this case, the Fisher--Bingham integral can be written by
the density function of $\chi$-distribution, and we have
$$
f(r) = \frac{2\pi^{d/2}}{\Gamma(d/2)} r^{d-1}e^{-r^2/2}.
$$
Table \ref{tab:comparison} shows the result for $f(1.0)$ by HGM 
and difference $\frac{2\pi^{d/2}}{\Gamma(d/2)} r^{d-1}e^{-r^2/2}-f(1.0)$
for each dimension.

\begin{table}[htbp]
\caption{Comparison to  $\chi$-distribution: $\Sigma=I$ and $\mu=0$.}
\label{tab:comparison}
\begin{center}
\begin{tabular}{ccc}
dim& hgm & exact$-$hgm \\
\hline
 3 & 7.621888 & 1.35e-06\\
 4 & 11.972435 & 7.09e-07\\
 5 & 15.963247 & 4.85e-07\\
 6 & 18.806257 & 3.70e-07\\
 7 & 20.060008 & 3.34e-07\\
 8 & 19.693866 & 3.13e-07\\
 9 & 18.005821 & 2.88e-07\\
10 & 15.467527 & 2.47e-07\\
\hline
\end{tabular}
\end{center}
\end{table}

As our fifth example we consider the case where 
$$
\Sigma=\mathrm{diag}\left(
\frac{1}{\sqrt{2}},\frac{1}{\sqrt{2}},
\frac{1}{\sqrt{4}},\frac{1}{\sqrt{4}},
\frac{1}{\sqrt{6}},\frac{1}{\sqrt{6}},
\cdots,
\frac{1}{\sqrt{2n}},\frac{1}{\sqrt{2n}}
\right),\,
\mu = \mathbf 0
\quad (d=2n).
$$
In this case, 
the ball probability \eqref{eq:cdf1} equals to 
$$
P\left(
 \frac{1}{2}X_1^2+\frac{1}{2}X_2^2
+\frac{1}{4}X_3^2+\frac{1}{4}X_4^2
+\frac{1}{6}X_5^2+\frac{1}{6}X_6^2
+\cdots
+\frac{1}{2n}X_{2n-1}^2+\frac{1}{2n}X_{2n}^2
< r^2
\right)
$$
where $X_1,\dots, X_{2n}$ are independent and identically distributed
with the standard normal distribution.
Since the distribution of $\frac{1}{2k}\left(X_{2k-1}^2+X_{2k}^2\right)$ is 
the exponential distribution with the rate parameter $k$,
the above probability is equal to $(1-e^{-r^2})^n$
\cite[p.21]{feller}.
The second column in Table \ref{tab:feller} shows 
the result of HGM for the ball probability at $r=1.0$.
The third column shows the difference between HGM and the exact value.

\begin{table}[htbp]
\caption{Comparison at specific parameters.}
\label{tab:feller}
\begin{center}
\begin{tabular}{ccc}
dim& hgm & exact$-$hgm \\
\hline
 6 & 0.252580 & 4.97e-09\\
 8 & 0.159661 & 2.54e-09\\
10 & 0.100925 & 1.61e-09\\
12 & 0.063797 & 1.03e-09\\
14 & 0.040327 & 8.16e-10\\
16 & 0.025492 & 7.07e-10\\
18 & 0.016114 & 3.04e-10\\
20 & 0.010186 & 2.37e-10\\
\hline
\end{tabular}
\end{center}
\end{table}

\section{Summary and discussion}
\label{sec:summary}

In this paper we applied HGM for computing
distribution function of a weighted sum of independent noncentral chi-square
random variables. We found that our method is numerically both
accurate and fast, after we implemented the following ideas. First,
during the application of Runge-Kutta method, we re-scaled the vector
$\mathbf F$ in \eqref{eq:vector-F} as needed to keep its elements within the
precision for floating point numbers.  Also we divided the interval for
integration into $(0,1]$ and $[1,\infty)$ and switched from $\mathbf F$
to $\mathbf Q$ in \eqref{eq:large-r-1} in view of the asymptotic values for
$\mathbf Q$.  Our experience in this paper 
shows that re-scaling of the standard monomials is important
in numerical implementation of HGM.

In our implementation, the numerical integration starts from a small $r=r_0 >0$ and
the integration proceeds to $r=\infty$.  On the other hand, we have asymptotic results
for large $r$ in Section \ref{sec:laplace}.  Then we might consider reversing the direction of
integration and start with initial values at very large $r$.   We may call the former the
``forward integration'' and the latter the ``backward integration''.
However we found that
the backward integration is not numerically stable.  Hence the asymptotic values
can not be used as initial values. 
In this paper we used the asymptotic values 
just for checking the accuracy HGM in the forward direction.

It is an interesting question, 
whether the asymptotic values can be used to adjust the
values of the forward integration.  We may look at the difference between $\mathbf F$ by forward HGM and
its asymptotic value for very large $r$ and use the difference to adjust $\mathbf F$ at intermediate values of $r$.
However it is not clear how this adjustment can be implemented.

\appendix
\section{A general form of Proposition \ref{prop:lap1}}

In Proposition \ref{prop:lap1} we assumed  $\x_1 > \x_2$.  
In this appendix we  state the following proposition
for the general case  $\x_1 = \dots = \x_m > \x_{m+1}$ without a proof.
For this case, the integrand for the Fisher-Bingham integral takes its maximum on the $(m-1)$-dimensional sphere $S^{m-1}(1)$, rather than on a finite number of points.  However by appropriate choice of coordinates
and by multiplication of the volume  ${\rm Vol}(S^{m-1}(1))$, 
the derivation of Proposition \ref{prop:lap-2} is basically the same as Proposition \ref{prop:lap1}.

\begin{proposition}
\label{prop:lap-2}
Assume that
\[
0 > \x_1 = \dots = \x_m > \x_{m+1} \ge \dots \ge \x_d.
\]
If $0=\y_1 = \dots = \y_m$, then as $r\rightarrow\infty$,
\begin{align*}
f(\x,\y,r) 
&=r^{m-1} S_{m-1} \exp\left( r^2 \x_1 -\sum_{i=m+1}^d \frac{\y_i^2}{4(\x_i - \x_1)} \right)
\frac{\pi^{(d-m)/2}}{\prod_{i=m}^d (\x_1-\x_i)^{1/2}} (1 + o(1)),  \\
\pd{\x_j} f(\x,\y,r) &= \frac{r^2}{m} f(\x,\y,r)(1+o(1)),  \quad j\le m, \\
\pd{\y_j} f(\x,\y,r) &= 0, \quad j\le m, \\
\pd{\y_j} f(\x,\y,r)
& = -\frac{\y_j}{2(\x_j-\x_1)} f(\x,\y,r)(1+o(1)), \quad j>m, \\
\partial_{\x_j} f(\x,\y,r)
&= \left(\frac{1}{2(\x_1-\x_j)} + \frac{\y_j^2}{4(\x_j-\x_1)^2}\right) f(\x,\y,r)(1+o(1)), \quad j>m.
\end{align*} 
If $(\y_1,\dots,\y_m) \neq (0,\dots,0)$, define 
$
\gamma = (\y_1^2 + \dots + \y_m^2)^{1/2}
$.
Then, as $r\rightarrow\infty$, 
\allowdisplaybreaks
\begin{align*}
f(\x,\y,r)&=
 \exp\left(r^2 \x_1 + r\gamma -\sum_{i=m+1}^d \frac{\y_i^2}{4(\x_i - \x_1)} \right)
\left(\frac{2r}{\gamma}\right)^{(m-1)/2}  
\frac{ \pi^{(d-1)/2}}{\prod_{i=m}^d (\x_1-\x_i)^{1/2}}(1+o(1)),\\
\partial_{\y_j} f(\x,\y,r)&= r \frac{\y_j}{\gamma} f(\x,\y,r)(1+o(1)), \quad \y_j \neq 0, \ j\le m, \\
\partial_{\x_j} f(\x,\y,r)&= r^2 \frac{\y_j^2}{\gamma^2} f(\x,\y,r)(1+o(1)), \quad \y_j\neq 0,\  j\le m,\\
\partial_{\y_j} f(\x,\y,r)&=0, \quad \y_j=0, \ j\le m,\\
\partial_{\x_j} f(\x,\y,r)&= \frac{r}{\gamma} f(\x,\y,r)(1+o(1)), \quad \y_j=0,\ j\le m, \\
\partial_{\y_j} f(\x,\y,r)& = -\frac{\y_j}{2(\x_j-\x_1)} f(\x,\y,r)(1+o(1)), 
\quad j>m, \\
\partial_{\x_j} f(\x,\y,r)
&= \left(\frac{1}{2(\x_1-\x_j)} + \frac{\y_j^2}{4(\x_j-\x_1)^2}\right) 
f(\x,\y,r)(1+o(1)), \quad j>m.
\end{align*}
\end{proposition}

\bibliographystyle{abbrv}
\bibliography{ball-probability}

\end{document}